\documentclass[12pt]{amsart}
\usepackage{latexsym,amssymb,enumerate, amsmath}

\newtheorem*{theoremA'}{Theorem A'}
\newtheorem*{theoremA"}{Theorem A"}

\usepackage{color}
\def\2zent#1{{\bf Z}_2 (#1)}

\begin{document}

\title{On  the  structure  of  some  locally  nilpotent  groups without  contranormal   subgroups}

\author[L.A. Kurdachenko]{Leonid A. Kurdachenko}
\address{L.A. Kurdachenko, Department of Algebra and Geometry \newline
School of Mathematics and Mechanics \newline
University of Dnipro\newline
Gagarin prospect 72,  Dnipro 10, \newline
49010 Ukraine}
\email{lkurdachenko@i.ua}

\author[P. Longobardi]{Patrizia Longobardi}
\address{P. Longobardi, Dipartimento di Matematica\newline
Universit\`a di Salerno\newline
via Giovanni Paolo II, 132, \newline
84084 Fisciano (Salerno), Italy}
\email{plongobardi@unisa.it}

\author[M. Maj]{Mercede Maj}
\address{M. Maj, Dipartimento di Matematica\newline
Universit\`a di Salerno\newline
via Giovanni Paolo II, 132, \newline
84084 Fisciano (Salerno), Italy}
\email{mmaj@unisa.it}


\begin{abstract} 
Following J.S. Rose, a subgroup $H$ of a group $G$ is said  contranormal  in $G$ if $G = H^G$. In a certain
sense, contranormal subgroups are antipodes to subnormal subgroups. It is well known that a finite group is
nilpotent if and only if it has no proper contranormal subgroups. We prove that a nilpotent-by-finite group with no proper contranormal subgroup is nilpotent.
There are locally nilpotent groups with a proper contranormal subgroup. We study the structure of hypercentral groups with a finite proper contranormal subgroup.
\end{abstract}

\maketitle

\smallskip

\begin{center}
  \emph{Dedicated to Professor Pavel Shumyatsky on his 60th birthday}
\end{center}

\bigskip

%%%%%%%%%%%%%%%%%%%%%%%%%%%%

{\bf 1. Introduction.}

\medskip

A subgroup $H$ of a group $G$  is called {\bf contranormal in $G$} if  $H^G = G$, where $H^G = \langle x^{-1}hx \ | 
\ h \in H,  \ x \in G\rangle$ is the normal closure of $H$ in $G$, the smallest normal subgroup of $G$ containing $H$. 
For example $G$ is contranormal in $G$, for any group $G$. The term  "contranormal subgroup" has been introduced by J.S. Rose in the paper  \cite{RJ1968}. Contranormal subgroups have been studied for example in the paper \cite{KOS2007}.
If $G$ is a group and $H$ is a contranormal subgroup of $G$, then every subgroup $K$ containing $H$ is contranormal in $G$. In particular, if $H$ and $L$ are contranormal subgroups of $G$, then the subgroup  $\langle H, L \rangle$  is also contranormal in $G$. However, the intersection of two contranormal subgroups is not always contranormal. For example, in the group  $A_4$ every Sylow  $3$-subgroup is contranormal, but the intersection of every two Sylow $3$-subgroups of $A_4$ is trivial, so that it is not contranormal. Notice also that if $M$  is a maximal subgroup of  $G$ which is not normal, then clearly $M$ is a contranormal subgroup of $G$.  Moreover, every subgroup of a finite group $G$ is a contranormal subgroup of a subnormal subgroup of $G$.  As we can see by the definition, contranormal subgroups are in a certain sense, antipode of normal and subnormal subgroups: a contranormal subgroup $H$ of a group $G$ is normal (respectively subnormal) if and only if $H = G$. It follows that groups, whose subgroups are subnormal (in particular, nilpotent group), do not contain proper contranormal subgroups. For finite groups the converse is true.

 {\it A finite group $G$ is nilpotent if and only if $G$  does not have proper contranormal subgroups. }

Indeed, suppose that there is a prime $p$ such that $G$ has a Sylow $p$-subgroup $P$ which is not normal in $G$. Then $N_G(P) \not= G$. Since $P$  is pronormal in $G$, $N_G(P)$ is abnormal in $G$ (\cite{RJ1967}, 1.6). But every abnormal subgroup is contranormal, and we obtain a contradiction, which shows that Sylow $q$-subgroups of $G$ are normal for each prime $q$. It follows that $G$ is nilpotent. 

There exist infinite non-nilpotent groups, whose subgroups are subnormal (it is possible to find examples of such groups in the survey \cite{CC2008}). Therefore the following question naturally appears:

 {\it When a locally nilpotent group without proper contranormal subgroups is nilpotent?}

We notice that there exist Chernikov  locally nilpotent groups having proper contranormal subgroups, as the following example shows.
Let $D$ be a divisible abelian $2$-group. Then $D$ has an automorphism $\varphi$ such that  $\varphi(d) = d^{-1}$ for each element  $d\in D$.  Define the semidirect product  $G = D \rtimes \langle b \rangle$  such that $d^b = \varphi(d) = d^{-1}$  for each element  $d \in D$.  Let $a$ be an arbitrary element of $D$. Since $D$ is divisible, there exists an element  $d\in D$ such that $d^2 = a$. We have $[b, d] = b^{-1} d^{-1} bd =  d^2 = a$. It follows that $[b, D] = D$. From $[b, D] \leq\langle  b \rangle^G$  and $\langle b \rangle \leq\langle b \rangle^G$  we obtain that  $\langle b \rangle^G  = \langle  b \rangle [b, D] = \langle b \rangle D = G$, so that the subgroup $\langle b \rangle$  is contranormal in $G$. We note that the group $G$  is not nilpotent, however the series 

$\langle 1 \rangle \leq\Omega_1(D) \leq \cdots \leq \Omega_n(D) \leq \Omega_{n+1}(D) \leq \cdots\leq D \leq G$

\noindent is central, so that $G$ is a hypercentral abelian-by-finite group. Besides, the contranormal subgroup  $\langle b \rangle$  is ascendant.
This group is abelian-by-finite, thus there exist hypercentral abelian-by-finite groups having  proper contranormal subgroups, and also finite contranormal subgroups.  
This example raises the following question:

{\it What can we say about locally nilpotent abelian-by-finite groups having  no proper contranormal subgroups?}

Our first result gives an answer to this question. In fact we have the following Theorem.
\vskip 0.4 true cm

{\bf Theorem A.}
{\it Let $G$ be a nilpotent-by-finite  group. If $G$ has no proper contranormal subgroups, then $G$ is nilpotent.}
\medskip
 
 Now  the question appears about the structure of locally nilpotent abelian-by-finite groups having proper contranormal subgroups.  We show here the following result.

\vskip 0.4 true cm
{\bf Proposition B.}
{\it Let $G$ be a locally nilpotent group and $A$ be a normal abelian subgroup of $G$ with $G/A$ finite.  Suppose that $G$ has a proper contranormal subgroup $C$, then 
$C = BK$ where $B \leq A$ is normal in $G$, $K$ is a finitely generated subgroup such that $G = AK$, and $A = B[K, A]$. In particular the factor group $G/B$ has the finite contranormal subgroup $KB/B$.}   
 \medskip

Therefore we naturally come to locally nilpotent abelian-by-finite groups having a finite contranormal subgroup. Our last result gives a description of hypercentral groups which include a finite contranormal subgroup.
\vskip 0.4 true cm

{\bf Theorem C.}
{\it Let $G$ be a hypercentral group. If $G$ contains a finite contranormal subgroup, then $G$ satisfies the following conditions:

$(i)$ $G = VC$, where $V$ is a normal divisible abelian subgroup and $C$ is a finite contranormal subgroup of $G$;

$(ii)$ $\Pi(G) = \Pi(C)$, in particular the set $\Pi(G)$ is finite;  

$(iii)$ $V$ has a family of $G$-invariant  $G$-quasifinite subgroups $\{D_\mu \ | \  \mu \in M\}$ such that  $V = \langle D_\mu \  | \ \mu \in M\}$;  

$(iv)$  $[D_\mu, C] = D_\mu$  for all $\mu \in M$, in particular, $[V, C] = V$. }

\medskip
Here an infinite normal abelian subgroup $A$ of a group $G$ is called {\bf $G$-quasifinite} if every proper $G$-invariant subgroup of $A$ is finite.

\bigskip

{\bf 2. Nilpotent-by-finite groups without proper contranormal subgroups}
\medskip

We start our investigation with this easy and very useful Lemma.
\medskip

{\bf Lemma 2.1.} \textit{Let $G$ be a group. Then:}

 \textit{$(i)$ If $C$ is a contranormal subgroup of $G$  and $H$ is a normal subgroup of $G$, then $CH/H$ is a contranormal subgroup of $G/H$. }
 
  \textit{$(ii)$ If $H$ is a normal subgroup of $G$ and $C$ is a subgroup of $G$ such that $H \leq C$ and $C/H$ is a contranormal subgroup of $G/H$, then $C$ is a contranormal subgroup of $G$. }

 \textit{$(iii)$ If $C$ is a contranormal subgroup of  $G$  and $D$ is a contranormal subgroup of  $C$, then $D$ is a contranormal subgroup of $G$.}
\begin{proof} 
These assertions are obvious. 
\end{proof}

Let $G$ be a nilpotent-by-finite group and assume that $G$ has no contranormal subgroups. In order to prove Theorem A, we first assume that $G$ is $p$-group, $p$ a prime.
Furthermore, we first suppose that $G$ is abelian-by-finite, thus there exists a normal abelian subgroup $A$ of $G$ of finite index in $G$.
We start  stating three easy Lemmas, well known in the literature. We add the proofs for the sake of completeness.
\medskip

{\bf Lemma 2.2.} \textit{Let $G$ be a $p$-group, $p$ a prime, and suppose that $G$ contains a normal bounded abelian subgroup $A$ such that $G/C_G(A)$ is finite. Then for some positive integer  $m$,
$A$ is contained in $\zeta_m(G)$, the $m-th$  term of the upper central series of $G$.}
\begin{proof}
Write $s$ the exponent of $A$ and $k = |G/C_G(A)|$. For each $a \in A$ we have $A \leq C_G(a)$ and $|G:C_G(a)| \leq k$. Thus $a$ has at most $k$ conjugates in $G$. Therefore $\langle a \rangle^G$  is an abelian group, 
of exponent $\leq s$, generated by at most $k$ elements. Thus $\langle a \rangle^G$ is a finite normal subgroup of order at most $s^k$. Write $m = s^k$. Since $G$ is a soluble $p$-group, then $G$ is locally nilpotent, hence $\langle a \rangle^G$
is contained in the $m-th$ term of the upper central series of $G$. That holds for each $a \in A$, therefore $A \leq \zeta_m(G).$
\end{proof}

{\bf Lemma 2.3.} \textit{Let $A$ be an abelian  $p$-group, $p$ a prime. If $A$ is not bounded, then $A$ contains a subgroup $B$ such that $A/B$ is a divisibile Chernikov group.}
\begin{proof}
Suppose first that $A$ is a direct product  of cyclic groups. Then since $A$ is not bounded, there exists a subgroup $C$ of $A$ such that $A/C = Dr_{n \in \mathbb{N}} \langle d_n \rangle$, 
where the element $d_n$ has order $p^n$. Consider the subgroup $B/C = \langle d_nd_{n+1}^{-p} \ | n \in \mathbb{N} \rangle$. Then by this choice the factor group $A/B$ is a Pr\"ufer $p$-group.

Suppose now that $A$ cannot be decomposed in a direct product of cyclic subgroups. Let $D$ be a basic subgroup of $A$ (see Theorem 32.3 of the book \cite{FL1970}). Then $D$ is the direct product of cyclic subgroups, therefore $D \not= A$. Moreover $A/D$ is a divisible group. Thus $A/D$ is direct product of Pr\"ufer $p$-groups and there exists a subgroup $B/D$ of $A/B$ such that $A/B$ is a Pr\"ufer $p$-group. 
\end{proof}

{\bf Lemma 2.4.} \textit{Let $\mathcal{H}$ be a class of groups closed under subgroups and under finite direct products.  
Let $G$ be a group containing a normal abelian subgroup $A$ such that $|G/C_G(A)|$ is finite. 
Suppose that $A$ contains a subgroup $B$ such that $A/B \in \mathcal{H}$, then $A$ contains a $G$-invariant subgroup $C$ such that $C \leq B$ and 
$A/C \in \mathcal{H}$.}
\begin{proof} 
For each element $g \in G$ the isomorphism $A/B^g \simeq A^g/B^g \simeq A/B$ shows that $A/B^g \in \mathcal{H}$. Since the subgroup $C_G(A)$ has finite index in $G$, the set $\{B^g \ | g \in G\}$ is finite. Write $\{B^g \ | \ g \in G\} = \{B_1, B_2, \dots , B_n\}$, and $C = B_1 \cap B_2\cdots \cap B_n$. Using Remak's theorem we obtain the embedding $A/C \lesssim A/B_1 \times A/B_2 \times\cdots \times A/B_n$. Since $ A/B_i \in \mathcal{H}$, for every $i \in \{1, \cdots, n\}$,  and $\mathcal{H}$ is closed under subgroups and finite direct products, it follows that $A/C \in \mathcal{H}$.
\end{proof}

Another general lemma we will use is the following:
\medskip

{\bf Lemma 2.5.} \textit{Let $G$ be a $p$-group, $p$ a prime, and suppose that $G$ contains a normal abelian subgroup $A$ such that $C_G(A)$ has finite index. Assume that $A$ contains a $G$-invariant  divisible Chernikov subgroup $D$. Then $A$ contains a $G$-invariant subgroup $S$ such that $A = SD$ and the intersection $S \cap D$ is finite.}

\begin{proof}
Since $D$ is divisible, it  has a complement in $A$, that is $A$ contains a subgroup $B$ such that $A = D \times B$. Then $A$ contains a $G$-invariant subgroup $C$ such that $(D \cap C)^n = \langle1\rangle$ and $A^n \leq DC$ where $n = |G/C_G(A)|$ (see, for example \cite{KOS2007}, Theoren 5.9). In particular, the intersection $D \cap C$ is finite. Then $DC/C \simeq  D/(D \cap C) \simeq D$. In particular, $DC/C$ is a divisible subgroup of $A/C$, therefore $A/C$ contains a subgroup $E/C$ such that $A/C = (DC/C) \times E/C$. Since the factor $A/DC$ is bounded, $E/C$ is bounded, moreover $(E/C)^n = \langle1\rangle$. Let $n = p^k$, then $E/C \leq \Omega_k(A/C)$. Put $S/C = \Omega_k(A/C)$, then the intersection $(S/C) \cap (DC/C)$ is finite and $A/C = (DC/C) (S/C)$. It follows that $A = DS$. Since $D \cap C$ and $(S/C)\cap (DC/C)$ are finite, then $S \cap D$ is finite. The Lemma is proved.

\end{proof}

Now assume that $G$ is a $p$-group, $p$ a prime, and that $G$ has no proper contranormal subgroups. Suppose that $G$ has a normal abelian subgroup $A$ of finite index in $G$. If $A$ is bounded, then
there exists a positive integer $m$ such that $A \leq \zeta_m(G)$, the $m-th$ term of the upper central series of $G$, by Lemma 2.2. Since $G/A$ is a finite $p$-group, $G/A$ is nilpotent. Therefore $G$ is nilpotent and we have the result of Theorem A in this case. Then we can suppose that $A$ is not bounded. Thus, by Lemma 2.3, there exists a subgroup $B$ of $A$ such that $A/B$ is a divisible Chernikov group. By Lemma 2.4 we can also suppose that $B$ is $G$-invariant. In this case we have.

\medskip
{\bf Lemma 2.6.} \textit{Let $G$ be a $p$-group, $p$ a prime, and suppose that $G$ contains a normal abelian subgroup $A$ of finite index. Assume that $A$ contains a $G$-invariant subgroup $C$ such that $A/C$ is a divisible Chernikov group.
If $G$ has no proper contranormal subgroups, then $[G, A] \leq C$.}

\begin{proof}
$A/C$ is a Chernikov group, thus $G/C$ satisfies the minimal condition on subgroups. Then there exists a series
\begin{center}

$C = C_1 \leq C_2 \leq \cdots \leq C_n = A$

\end{center}

\noindent of $G$-invariant subgroups such that the factors $C_{j+1}/C_j$ are $G$-quasifinite, $j \in \{1, \dots, n\}$.
Consider the factor $A/C_{n-1}$.  The subgroup $[G/C_{n-1}, A/C_{n -1}]$ is $G$-invariant, then either $[G/C_{n -1}, A/C_{n -1}] = A/C_{n-1}$, or $[G/C_{n-1}, A/C_{n-1}]$ is finite. Assume that $[G/C_{n-1}, A/C_{n -1}] = A/C_{n-1}$. Choose a finite subgroup $K/C_{n-1}$ such that $G/C_{n-1} = (A/C_{n-1})(K/C_{n-1})$. Then $[G/C_{n-1}, A/C_{n-1}] = [K/C_{n-1}, A/C_{n -1}]$. Then the inclusion $A/C_{n-1} = [K/C_{n-1}, A/C_{n-1}] \leq  (K/C_{n-1})^{G/C_{n-1}}$  implies that 
$(K/C_{n-1})^{G/C_{n-1}} = (A/C_{n-1})(K/C_{n-1}) = G/C_{n-1}$. This means that the subgroup $K/C_{n-1}$ is contranormal in $G/C_{n -1}$. By Lemma 2.1, the subgroup $K$ is contranormal in $G$, and we obtain a contradiction. 
This contradiction shows that $[G/C_{n-1}, A/C_{n-1}]$ is finite.  In this case the factor group $G/C_{n-1}$ is nilpotent. It follows that the center of $G/C_{n-1}$ contains $A/C_{n-1}$ (see, for example, \cite{DKS2017}, Proposition 3.2.11). 
Hence $[G, A] \leq C_{n -1}$.

\noindent Suppose that we have already proved that $[G, A] \leq C_2$. Since the subgroup $A/C$ is divisible and Chernikov, $A/C$ contains a $G$-invariant divisible subgroup $D/C$ such that $A/C 
= (C_2/C)(D/C)$ and the intersection $(C_2/C) \cap (D/C)$ is finite (see, for example, \cite{KOS2007}, Corollary 5.11). Then the factor $A/D$ is divisible Chernikov and $G$-quasifinite. Using the result of the previous paragraph, we obtain that 
$[G, A] \leq D$. Thus we have $[G/C, A/C] \leq C_2/C$ and $[G/C, A/C] \leq D/C$, therefore $[G/C, A/C] \leq (C_2/C) \cap (D/C)$. Since the last intersection is finite, the factor group $G/C$ is nilpotent. 
It follows that the center of $G/C$ contains $A/C$ (see, for example, \cite{DKS2017}, Proposition 3.2.11). Hence $[G, A] \leq C$, and the Lemma is proved.
\end{proof}

From Lemma 2.6 we have the following lemma:
\medskip

{\bf Lemma 2.7.} \textit{Let $G$ be a $p$-group, $p$ a prime, and suppose that $G$ contains a normal abelian subgroup $A$ of finite index. 
If $G$ has no contranormal subgroups, then $[G, A]$ is bounded.}
\begin{proof}
If $A$ is bounded, we have the result. Therefore we suppose that $A$ is not bounded. Then Lemma 2.3 shows that $A$ contains a subgroup $B$ such that $A/B$ is a divisible Chernikov group. Then Lemma 2.6 implies that $[G, A] \not= A$, moreover $A/[G, A]$ is not bounded. Suppose that the  subgroup $D = [G, A]$ is not bounded. Using again Lemma 2.3 we obtain that $D$ contains a subgroup $C$ such that $D/C$ is a divisible Chernikov group. Then, by Lemma 2.4, there exists a
 $G$-invariant subgroup $E$  such that $D/E$ is a Chernikov group. Then $D$ contains a $G$-invariant subgroup $H$ such that $E \leq H$, $H/E$ is finite and $D/H$ is a divisible Chernikov group. Therefore without loss of generality we may suppose that $D/E$ is a divisible Chernikov group.
We have $[G/E, A/E] = [G, A]E/E = DE/E = D/E$. Therefore $[G/E, A/E]$ is a divisible Chernikov group. By Lemma 2.5, $A/E$ contains a $G$-invariant subgroup $S/E$ such that $A/E = (D/E)(S/E)$ and the intersection $(D/E) \cap (S/E)$ is finite. It follows that $A/S \simeq (A/E)/(S/E) = (D/E)(S/E)/(S/E) \simeq (D/E)/((D/E) \cap (S/E)) \simeq D/E$ is a divisible Chernikov group. Furthermore, $A/S = (DS)/S = [G, A]S/S = [G/S, A/S]$. Now, by Lemma 2.6, $[G/E, A/E]\leq S/E$, since $A/S$ is a divisible Chernikov. Then $[G, A] \leq S$  and we obtain the contradiction $A = S$. This contradiction proves that the subgroup $[G, A]$ is bounded.

\end{proof}

Now we can prove the result of Theorem A, if $G$ is an abelian-by-finite $p$-group, $p$ a prime. 

\medskip

{\bf Corollary 2.8.} \textit{Let $G$ be a $p$-group, $p$ a prime, and suppose that $G$ contains a normal abelian subgroup $A$ of finite index. 
If $G$ has no proper contranormal subgroups, then $G$ is nilpotent.}
\begin{proof}
By Lemma 2.7, $[G, A]$ is bounded. Then, by Lemma 2.2, there exists a positive integer $t$ such that $[G, A] \leq \zeta_t(G)$. Then $A \leq \zeta_{t+1}(G)$, and $G$ is nilpotent since $G/A$ is a finite $p$-group.
\end{proof}

Next step is to  prove the result of Theorem A for every locally nilpotent abelian-by-finite group.
\medskip

{\bf Corollary 2.9.} \textit{Let $G$ be a locally nilpotent group, and suppose that $G$ contains a normal abelian subgroup $A$ of finite index. 
If $G$ has no proper contranormal subgroups, then $G$ is nilpotent.}
\begin{proof}
First, suppose that $G$ is periodic. Let $\pi = \Pi(G/A)$ and $\sigma = \Pi(G) \setminus \pi$, then the set $\pi$ is finite and we have $G = Dr_{p \in \pi} G_p \times Dr_{p \in \sigma} G_p$, where $G_p$ is a Sylow $p$-subgroup of $G$ for all $p \in \Pi(G)$. The isomorphism $G_p \simeq G/Dr_{q \in \Pi(G), q \not= p} G_q$ and Lemma 2.1 show that $G_p$ has no proper contranormal subgroups for every $p \in \pi$. Using Corollary 2.8 we obtain that $G_p$ is nilpotent for each $p \in \pi$. The finiteness of the set $\pi$ implies that $Dr_{ p \in \pi} G_p$ is nilpotent. Obviously the subgroup $G_p$ is abelian for every $p \in \sigma$, hence $Dr_{ p \in \sigma} G_p$ is abelian. Therefore $G$ is nilpotent.
Now suppose that $G$ is non-periodic. Then the set $Tor(G)$ of all elements of $G$ having finite order, is a characteristic subgroup of $G$ and the factor  group $S = G/Tor(G)$ is torsion-free. On the other hand, $S$ is abelian-by-finite. then $S$ is a locally nilpotent torsion-free abelian-by-finite group, and then it is abelian (see, for example, \cite{DKS2017}, Corollary 1.2.8). Choose in the abelian subgroup $A$ a  maximal $\mathbb{Z}$-independent subset $M$ and let $C$ be the subgroup of $A$ generated by $M$. Then $A/C$ is a periodic group. By Lemma 2.4 there exists a $G$-invariant subgroup $E \leq C$ such that $A/E$ is periodic. Obviously $E$ is torsion-free. Then $E \cap Tor(G) = \langle 1\rangle$. Using Remak's theorem, we obtain an embedding $G \lesssim G/E \times G/Tor(G)$. By Lemma 2.1 $G/E$ does not include proper contranormal subgroups. Then $G/E$ is nilpotent by Corollary 2.8, moreover  $G/Tor(G)$ is abelian, therefore $G$ is nilpotent and we have the result.
\end{proof}

Now we extend Corollary 2.9 to any abelian-by-finite group. We start with the following two results.

\medskip

{\bf Lemma 2.10.} \textit{Let $G$ be a group and suppose that $G$ contains a normal abelian $p$-subgroup $A$ of finite index, where $p$ is a prime. 
If $G$ has no proper contranormal subgroups, then $G$ is nilpotent.}
\begin{proof}
By Lemma 2.1 the factor group $G/A$ does not contain proper contranormal subgroups. Being finite, $G/A$ is nilpotent. Then $G/A = P/A \times S/A$,  where $P/A$ is a $p$-group and $S/A$ is a $p^{\prime}$-group. We have $A = C_A(S) \times [S, A]$ (see, for example, \cite{BKOP2008}, Proposition 2.12).
Suppose that the subgroup $[S, A]$ is not trivial. Since the subgroup $S$ is normal in $G$, then both subgroups $C = C_A(S)$ and $[S, A]$ are $G$-invariant. Moreover, we have $A/C = C[S, A]/C = [S/C, A/C]$.
If the abelian $p$- group $A/C$ is bounded, then it is the direct product of cyclic subgroups. In particular, $A/C$ contains a proper subgroup having finite index. Then, by Lemma 2.4, $A/C$ contains a proper $G$-invariant subgroup $B/C$, having finite index. By Lemma 2.1 the factor group $G/C$ does not contain proper contranormal subgroups. Being finite, this factor group must be nilpotent. But in this case $[A/C, S/C] = \langle 1\rangle$, and we obtain a contradiction.
If the abelian $p$-group $A/C$ is not bounded, then by Lemma 2.3, $A/C$ contains a subgroup $D/C$ such that $A/D$ is a divisible Chernikov group. By Lemma 2.4, $A/C$ contains a proper $G$-invariant subgroup $E/C$ such that $A/E$ is Chernikov. By Lemma 2.1 the factor group $G/E$ does not contain proper contranormal subgroups. Being Chernikov, this factor group must be nilpotent  (\cite{KS2003}, Lemma 4.9). But in this case $[A/C, S/C] = \langle 1\rangle$, and we again obtain a contradiction.
This contradiction proves that $A = C_A(S)$. It follows that $S = A \times V$ where $V$ is a finite $p^{\prime}$-subgroup. Moreover, $V$ is a Sylow $p^{\prime}$-subgroup of $S$, so that $V$ is normal in $G$. By Lemma 2.1 the factor group $G/V$ does not contain proper contranormal subgroups. This factor group is an abelian-by-finite $p$-group, then it is nilpotent, by Corollary 2.8. The equality $A \cap V = \langle1\rangle$ and Remak's theorem imply the embedding $G \lesssim  G/A \times G/V$, which implies that $G$ is nilpotent.
\end{proof}

Let $G$ be a group and $A$ be a normal subgroup of $G$. We put  $\gamma_(G, A) = A$, $\gamma_2(G, A) = [G, A]$, and, recursively, $\gamma_{\alpha+1}(G, A) = [G, \gamma_\alpha(G, A)]$, for all ordinals $\alpha$, moreover, if $\lambda$
is a limit ordinal,  we write $\gamma_\lambda(G, A) = \bigcap_{\mu < \lambda} \gamma_\mu(G, A)$ 

\medskip
{\bf  Lemma 2.11.} \textit{Let $G$ be a group and suppose that $G$ contains a normal abelian torsion-free subgroup $A$ of finite index.
If $G$ has no proper contranormal subgroups, then $G$ is nilpotent.}
\begin{proof}
Let $M$ be a finite subset of $A$ and write $B = \langle M \rangle^G$. Since $G/A$ is finite, the subgroup $B$ is finitely generated. Being torsion-free, it is free abelian. Moreover, $B$ is $G$-invariant. Put $T/B = Tor(A/B)$, then the subgroup $T$ has finite $0$-rank and it is $G$-invariant. Let $r_0(T) = n$, then $T/B$ has special rank at most $n$. Let $p$ be an arbitrary prime and consider the factor $A/B^p$. Let $S_p/B^p$ be the Sylow $p$-subgroup of $A/B_p$, then $S_p/B^p$ is  a Chernikov group of special rank at most $n$. We have the direct decomposition $A/B^p = Sp/B^p \times Cp/B^p$ (see, for example \cite{FL1970}, Theorems 21.2 and 27.5). Thus $A/C^p$ is a Cernikov $p$-group of special rank at most $n$. By Lemma 2.4  there exists a $G$-invariant subgroup $D_p$, $D_p \leq C_p$ such that $A/D_p$ is a Chernikov $p$-group, it is $D_p =  \bigcap_{g \in G} C_p^g$, thus $A/D_p$ has special rank at most $kn$ where $k = |G/A|$. The inclusion $D_p \leq C_p$ implies that 
$B \cap D_p = B^p$. It follows that $(BD_p)/D_p \simeq B/(B \cap D_p) = B/B^p$, in particular $(BD_p)/D_p$ is an elementary abelian $p$-group, having finite order less or equal to $p^n$. The factor-group $G/D_p$ is periodic, therefore, by Corollary 2.9, $G/D_p$ is nilpotent. Then $(BD_p)/D_p \leq \gamma_{n}(G)$, the $n-th$ term of the lower central series of $G$. It follows that $\gamma_{n+1}(G, B) \leq D_p$. On the other hand, since $B$ is normal in $G$, $\gamma_{n+1}(G, B) \leq B$, so that $\gamma_{n+1}(G, B) \leq D_p \cap B = B^p$. The last inclusion is true for each prime $p$, therefore $\gamma_{n+1}(G, B) \leq\bigcap_{p \in  P} B^p$, where $P$ is the set of all primes.. Since $B$ is a free abelian subgroup, $\bigcap_{p \in P} B^p = \langle 1 \rangle$, thus $\gamma_{n+1}(G, B) = \langle1\rangle$. It follows that $B \leq \gamma_{n}(G )$. That holds for every finitely generated subgroup $B$ of $A$, therefore $A$ is contained in the hypercenter of $G$. By Lemma 2.1 the factor  group $G/A$ does not contain proper contranormal subgroups. Being finite, $G/A$ is nilpotent. Then $G$ is hypercentral. In particular, $G$ is locally nilpotent, and, by Lemma 2.9, $G$ is nilpotent.
\end{proof}

{\bf  Corollary 2.12.} \textit{Let $G$ be an abelian-by-finite group. If $G$ has no proper contranormal subgroups, then $G$ is nilpotent.}
\begin{proof}
Let $A$ be a normal abelian subgroup of $G$ such that the factor  group $G/A$ is finite. First suppose that $G$ is periodic. Let $\pi = \Pi(G/A)$ and $\sigma = \Pi(A) \setminus \pi$, then the set $\pi$ is finite and we have 
$A = Dr_{ p \in \pi} A_p \times Dr_{p \in  \sigma} A_p$, where $A_p$ is the Sylow $p$-subgroup of $A$ for all $p \in \Pi(A)$. Put $B_p = Dr_{q \in \Pi(A), q \not= p} A_q$, then the subgroup $B_p$ is $G$-invariant, $A/B_p \simeq A_p$
 and by Lemma 2.1 $G/B_p$ does not contain proper contranormal subgroups for every $p \in \Pi(A)$. By Lemma 2.10 $G/B_p$ is nilpotent for each $p \in \Pi(A)$. In particular, if $p \in \sigma$, then $G/B_p$ is abelian. 
 Since $\bigcap_{ p \in \Pi(A) }B_p = \langle 1\rangle$, by Remak's theorem, we obtain an embedding $G \lesssim Dr_{p \in \pi} G/B_p \times Cr_{ p \in  \sigma} G/B_p$. Since the set $\pi$ is finite $Dr_{p \in \pi} G/B_p$ is nilpotent. 
 Since $G/B_p$ is abelian for all $p \in \sigma$, then $Cr_{p \in \sigma} G/B_p$ is abelian. Therefore $G$ is nilpotent. Now suppose that $G$ is not periodic. Since $G$ is not periodic, $A$ is also not periodic. write $T = Tora(A)$. 
 Then $A \not =  T$. Obviously the subgroup $T$ is $G$-invariant and $A/T$ is torsion-free.  Lemma 2.1 shows that $G/T$ does not contain proper contranormal subgroups. Hence the factor group $G/T$ is nilpotent, by Lemma 2.11.
Choose in the abelian subgroup $A$ a maximal $\mathbb{Z}$-independent subset $M$ and let $C = \langle M \rangle$. Then $A/C$ is a periodic group. By Lemma 2.4 there exists a $G$- invariant  subgroup $E$ such that $E \leq C$ and $A/E$ is a periodic group. The inclusion $E \leq C$ implies that $E$ is torsion-free. Thus $E \cap T = \langle 1\rangle$. By  RemakÕs theorem, we obtain an embedding $G \leq G/E \times G/T$. Lemma 2.1 shows that $G/E$ does not contain proper contranormal subgroups. Being periodic, $G/E$ is nilpotent, we know that $G/T$ is nilpotent, hence $G$ is nilpotent, as required.
 
\end{proof}

Now we can prove Theorem A.
\medskip

{\it Proof of Theorem A}.
 Let $K$ be a nilpotent normal subgroup of $G$ such that $G/K$ is finite. Write $D = [K, K]$. Lemma 2.1 implies that the factor group $G/D$ does not contain proper contranormal subgroups. Moreover, $G/D$ is abelian-by- finite. Then Corollary 2.12 implies that $G/D$ is nilpotent. Using now Theorem 7 of paper \cite{HP1958}, we obtain that $G$ is nilpotent, as required.
\qed

\medskip
{\bf 3. Locally nilpotent abelian-by-finite groups with a finite contranormal subgroup}

\medskip

We start this section by proving Proposition B.

\medskip
{\it Proof of Proposition B.}
Suppose that $AC \not = G$. Then Lemma 2.1 implies that $CA/A$ is a proper contranormal subgroup of the finite nilpotent group $G/A$. But a nilpotent group does not contain a proper contranormal subgroups. Hence $AC = G$. Choose in $C$ a finitely generated subgroup $K$ such that $AK = G$, then $C = BK$ where $B = C \cap A$.
Since $A$ is normal in $G$, $[K, B] \leq A$. On the other hand, $[K, B] \leq C$, so that $[K, B] \leq C \cap A = B$. Therefore, the subgroup $B$ is $K$-invariant. $B$ is also $A$- invariant, since $A$ is abelian, thus from $G = AK$ we get that $B$ is $G$-invariant. The intersection $K \cap A$ is normal in $G$. Considering the factor group $G/(K \cap A)$,  without loss of generality we may assume that $K \cap A$ is trivial. Then the subgroup $K$ is finite.
From $G = AK$, with $A$ normal in $G$,  it follows $[K, A]$ normal in $G$ and $[G, G] = [K, A][K, K] \leq K[A,K]$ Thus $G/(K[K, A])$ is abelian. By Lemma 2.1 $C[K, A]/(K[K, A])$ is contranormal in $G/(K[K, A])$. It follows that $C[K, A]/(K[K, A]) = G/(K[K, A])$. Therefore we have $G = C[K, A] = BK[K, A] = B[K, A] \rtimes K$. In particular, we obtain that $A = B[K, A]$.
The subgroup $B$ is normal in $G$. Then we obtain that $G/B = A/B \rtimes KB/B = [K, A]B/B \rtimes KB/B = [KB/B, A/B] \rtimes KB/B.$
It follows $G/B = (KB/B)^{G/B}$, hence $KB/B$ is contranormal in $G/B$.
\qed
\medskip

We start our investigation assuming that $G$ is a $p$-group, $p$ a prime.

\medskip

{\bf Proposition 3.1.} \textit {Let $G$ be an abelian-by-finite $p$-group, $p$ a prime. If $G$ contains a finite contranormal subgroup, then $G$ satisfies the following conditions:}

 \textit{$(i)$  G = VC  where  V  is a normal divisible abelian subgroup and  $C$  is a finite contranormal subgroup of  $G$;}

 \textit{$(ii)$  $V$  has a family of $G$-invariant $G$-quasifinite subgroups $\{ D_\mu \ | \ \mu \in M \}$ such that $V = \langle D\mu \  | \  \mu \in M \rangle$;}
 
 \textit{$(iii)$  $[D\mu, C] = D\mu$ for all $\mu \in M$, in particular, $[V, C] = V$. }

\begin{proof}   Let  $A$  be a normal abelian subgroup of $G$ having finite index and let  $C$  be a finite contranormal subgroup of $G$. By Lemma 2.1 $CA/A$  is contranormal in $G/A$. Since $G/A$ is a finite $p$-group, it is nilpotent. The fact that a nilpotent group does not include proper contranormal subgroups implies that  $CA/A = G/A$  or $G = CA$. 
 If $A = A^p$, then $A$ is divisible and $(i)$ holds. Suppose that $B = A^p \not= A$. Then $B$ is normal in $G$ and $G/B$ is an extension of an elementary abelian $p$-subgroup by a finite $p$-group. Such groups are nilpotent (\cite{BG1959}). On the other hand, Lemma 2.1 shows that $CB/B$ is a contranormal subgroup of $G/B$. The fact that a nilpotent group does not include a proper contranormal subgroup implies that $CB/B = G/B$. It follows that  $A/B$  is finite. The finiteness of  $A/A^p$   implies that  $A = F \times V$  where  $V$  is a divisible subgroup and $F$  is a finite subgroup (see, for example  \cite{KL1986}, Lemma 3). Clearly the subgroup $V$ is $G$- invariant. Being a finite $p$-group, the factor group $G/V$  is nilpotent. As above it follows that $CV/V = G/V$  or  $G = VC$, ad again $(i)$ holds.
 Now suppose $G = VC$, where $V$ is divisible, abelian and normal in $G$.
Since $V$  is an abelian divisible $p$-subgroup we have  $V = \times_{ \lambda \in \Lambda}  P_\lambda$, where $P_\lambda$ is a Pr\"ufer $p$-subgroup for all $\lambda \in \Lambda$ (see, for example  \cite{FL1970}, Theorem 23.1). 
Let $Q_1$ be a Pr\"ufer $p$-subgroup of  $V$. Since $G/V$ is finite, $Q_1$ has only finitely many conjugates, so that $Y = Q_1^G$ is a divisible Chernikov subgroup. Since $Y$ satisfies the minimal condition, $Y$ includes an infinite $G$-invariant subgroup $D_1$ which  is $G$-quasifinite. If $D_1^p \not= D_1$, then $D_1^p$ is finite since $D_1$ is quasi finite, and $D_1/D_1^p$ is finite since it is an elementary abelian $p$-group with the minimal condition, hence $D_1$ is finite, a contradiction. Therefore $D_1^p = D_1$ and $D_1$ is divisible.  Thus $V = D_1R$  for some  subgroup $R$ such that $R$ is $G$-invariant, the intersection $D_1 \cap R$ is finite and $(D_1 \cap R)^{|C|} = \langle 1\rangle$  (see, for example \cite{KOS2007}, Corollary 5.11]. Put $|C|=p^n$, then $D_1\cap R \leq \Omega_n(V)$.
 It is not hard to prove that the subgroup $[D_1, C]$ is $G$-invariant. If we suppose that $[D_1, C]$ is a proper subgroup of $D_1$, then the fact that $D_1$ is $G$-quasifinite implies that $[D_1, C]$ must be finite. Then $D_1C$ is a finite-by-abelian $p$-group, so that $D_1C$ is nilpotent. Being Chernikov, $D_1C$ is central-by-finite (see, for example \cite{DKS2017}, Corollary 3.2.10). It follows that $D_1\leq  \zeta(D_1C)$. 
 Consider the factor group $G/R$. We have $V/R = D_1R/R\simeq D_1/(D_1 \cap R)$. The equality $[D_1, C] = \langle1\rangle$  implies that  $[V/R, CR/R] = [D_1R/R, CR/R] = [D_1, C]R/R = \langle 1\rangle.$ It follows that $V/R \leq \zeta(G/R)$. But in this case $(CR/R)^{G/R} = CR/R$, and we obtain a contradiction with Lemma 2.1. This contradiction shows that $[D_1, C]  = D_1$. 
 Choose in the subgroup $R$ a  Pr\"ufer $p$-subgroup $Q_2$. Again $Q_2$ has only finitely many conjugates, so that $Q_2^G$ is a divisible Chernikov subgroup. As above $Q_2^G$ includes an infinite $G$-invariant subgroup $D_2$, which  is $G$-quasifinite. Arguing as before it is possible to prove that $D_2$ is divisible. Then, by Corollary 5.11 of \cite{KOS2007}, $R = D_2R_1$ for some subgroup $R_1$ such that $R_1$ is $G$-invariant and the intersection $D_2 \cap R_1$ is finite, moreover $D_2 \cap R_1 \leq \Omega_n(V)$.
Using the above arguments, we obtain that $[D_2, C] = D_2$. 
Put $L_1 = \Omega_n(D_1)$, then $D_1/L_1 \cap RL_1/L_1 = \langle 1\rangle$  and $L_1 \leq \Omega_n(V)$. Similarly, put $L_2 = \Omega_n(D_2)$, then $D_2/L_2 \cap R_1L_2/L_2 = \langle 1\rangle$ and $L_2 \leq \Omega_n(V)$. Repeating these arguments and using transfinite induction, we obtain that the subgroup $V$ has a family of $G$-invariant $G$-quasifinite subgroups  $\{ D_\mu \ |\  \mu \in M \}$ such that  $V = \langle D_\mu \ | \mu \in M \rangle$, $[D_\mu, C] = D_\mu$ for all $\mu \in M$, as required. Moreover we have $V/\Omega_n(V) = \times_{\mu \in M} D_\mu \Omega_n(V)/\Omega_n(V).$
\end{proof}
\medskip

Now we can prove 

\medskip

{\bf Corollary 3.2.}
{\it Let $G$ be a periodic locally nilpotent abelian-by-finite  group. If $G$ contains a finite contranormal subgroup, then $G$ satisfies the following conditions:}

{\it $(i)$ $G = V C$  where $V$ is a normal divisible abelian subgroup and $C$ is a finite contranormal subgroup of $G$;}

{\it $(ii)$ $\Pi(G) = \Pi(C)$, in particular the set $\Pi(G)$ is finite; }

{\it $(iii)$ $V$ has a family of $G$-invariant  $G$-quasifinite subgroups $\{D_\mu \ | \  \mu \in M\}$ such that  $V = \langle D_\mu \  | \ \mu \in M\}$;  }

{\it $(iv)$  $[D_\mu, C] = D_\mu$  for all $\mu \in M$, in particular, $[V, C] = V$. }

\begin{proof} Let $A$ be a normal abelian subgroup of $G$ having finite index, and let $C$ be a finite contranormal subgroup of $G$. Then, arguing as above we have $G = CA$. 
Suppose that $\Pi(G) \not= \Pi(C)$ and choose a prime $q \in \Pi(G) \setminus \Pi(C)$.  The equality $G = AC$ implies that $A$ contains a Sylow $q$-subgroup $Q$ of $G$. We have $A = Q \times R$ where $R$ is a Sylow $q^{\prime}$- subgroup  of $A$. Then $G/R = QR/R \times CR/R$, which shows that $CR/R$ cannot be a contranormal subgroup of $G/R$. Thus we obtain a contradiction with Lemma 2.1. This contradiction proves that  $\Pi(G) = \Pi(C)$. We have $G = \times_{ p \in \Pi(G)} S_p$ where $S_p$ is a Sylow $p$-subgroup of $G$. The isomorphism  $S_p \simeq G/(\times_{ q \in \Pi(G), q \not= p} S_p)$ and an application of Proposition 3.1 prove the result. 
\end{proof}

\medskip

Recall that a group $G$ is called $\mathcal{F}$-perfect, if $G$ does not contain a proper subgroup of finite index. In every group the subgroup $\mathcal{F}(G)$, generated by all $\mathcal{F}$-perfect subgroups, is $\mathcal{F}$-perfect. It is the greatest $\mathcal{F}$-perfect subgroup of $G$. Clearly $\mathcal{F}(G)$ is a characteristic subgroup of $G$, and the factor group $G/\mathcal{F}(G)$ does not contain $\mathcal{F}$-perfect subgroups. The subgroup $\mathcal{F}(G)$ is called the $\mathcal{F}$-perfect part of $G$.
Let $\mathcal{X}$ be a class of groups. If $G$ is a group, then we denote by $G_{\mathcal{X}}$ the intersection of all normal subgroups $H$ of $G$ such that $G/H \in \mathcal{ X}$. The subgroup $G_{\mathcal{X}}$ is called the $\mathcal{X}$-residual of the group $G$. If $\mathcal{X} = \mathcal{F}$ is the class of all finite groups, then $G_{\mathcal{F}}$ is called the finite residual of $G$.

\medskip
{\bf Lemma 3.3.} {\it Let $G$ be a locally nilpotent periodic group. If $G$ contains a finite contranormal subgroup, then the $\mathcal{F}$-perfect part of $G$ has finite index.}

\begin{proof} If $G$ does not contain proper subgroups of finite index, then $G$ is $\mathcal{F}$-perfect and the result is proved. Therefore we suppose that $G$ contains proper subgroups of finite index.
Let $S$ be a finite contranormal subgroup of $G$. Then $S$ is nilpotent. Let $k$ be the nilpotency class of $S$. If $H$ is a normal subgroup of $G$ such that $G/H$ is finite, then Lemma 2.1 shows that $SH/H$ is a contranormal subgroup of $G/H$. On the other hand, $G/H$ is nilpotent, and a nilpotent group does not contain proper contranormal subgroups. It follows that $SH/H = G/H$. In particular, $G/H$ has nilpotency class at most $k$.
Let $\mathcal{S}$ be the family of all normal subgroups of $G$ having finite index, and let $L =  \bigcap_{ H \in \mathcal{ S}} H$. By Remak's theorem there is an embedding $G/L \lesssim  Cr _{H \in \mathcal {S} }G/H$. Since $G/H$ has nilpotency class at most $k$ for every $H \in \mathcal{S}$, this implies that $G/L$ is a nilpotent group. It follows that $G/L$ does not contain proper contranormal subgroups and we obtain the equality $G/L = SL/L$. This means that $G/L$ is finite.
If we suppose that $L$ contains a proper subgroup $K$ having finite index in $L$, then $K$ has finite index in $G$. Then $D = Core_G(K)$ is normal in $G$ and has finite index in $G$. Then $D \in \mathcal {S}$, and therefore $L \leq D$, a contradiction. This contradiction proves that $L$ is $\mathcal{F}$-perfect and $L$ coincides with the $\mathcal{F}$-perfect part of $G$.
\end{proof}

\medskip

{\bf Corollary 3.4.} {\it Let $G$ be a hypercentral periodic group. If $G$ contains a finite contranormal subgroup, then $G$ is abelian-by-finite.}
\begin{proof} Let $L$ be the $\mathcal{F}$-perfect part of $G$. Lemma 3.3 implies that $L$ has finite index in $G$. The result follows since a periodic hypercentral $\mathcal{F}$-perfect group is abelian (see \cite{CSN1980}, Chapter 2, n.  2, Theorem 2.2).
\end{proof}

\medskip 

{\bf Lemma 3.5.} {\it Let $G$ be a locally nilpotent group. If $G$ is not periodic, then $G$ does not contain finite contranormal subgroups. }

\begin{proof}  Suppose the contrary, and let $S$ be a finite contranormal subgroup of $G$.  Since $G$ is locally nilpotent, the set $Tor(G)$ of all elements of $G$ having finite order is 
a characteristic subgroup of $G$. Since $G$ is not periodic, $G \not= Tor(G)$. Then the inclusion $S \leq Tor(G)$ implies that $S^G \not= G$ and we obtain a contradiction which proves the result.

\end{proof}

\medskip

Now we can prove Theorem C.

\medskip

{\it Proof of Theorem C.}
Lemma 3.5 implies that a group $G$ must be periodic. By Corollary 3.4 $G$ is abelian-by-finite, and the result follows from Corollary 3.2.
\qed

\bigskip

\section*{Acknowledgements}

This work was supported by the ``National Group for Algebraic and Geometric 
Structures, and their Applications" (GNSAGA - INdAM), Italy.
 
 \medskip

The first author is grateful to the Department of Mathematics of the University 
of Salerno for its hospitality and support, while this investigation was carried out.
\medskip

After writing this work we noticed that Theorem A was also proved by B.A.F. Wehrfritz in the paper "Groups with no proper contranormal subgroups", Publ. Mat. 64 (2020) 183-194.

\bigskip
\end{document}